\documentclass{article}

\input amssym.def

\makeatletter
\@addtoreset{equation}{section}
\makeatother

\newtheorem{prop}{Proposition}
\newtheorem{exam}{Example}

\def\be{\begin{equation}}
\def\ee{\end{equation}}

\begin{document}

\begin{titlepage}

\title{On pentagon, ten-term, and tetrahedron relations}

\author{R.M. Kashaev\thanks{On leave of absence from
St. Petersburg Branch of the Steklov Mathematical Institute}
\thanks{Supported by MAE-MICECO-CNRS Fellowship}\\ \\
Laboratoire de Physique Th\'eorique 
{\sc enslapp}\thanks{URA 14-36 du CNRS,
associ\'ee \`a l'E.N.S. de Lyon, 
et \`a l'Universit\`e de Savoie}
\\
ENSLyon, 46 All\'ee d'Italie,\\
69007 Lyon, FRANCE\\
\\
S.M. Sergeev\thanks{Partially supported by INTAS 
Grant 93-2492 and RFFR Grant
95-01-00249}\\ \\
Scientific Center of Institute of Nuclear Physics SB RAS\\
Protvino, Moscow Region 142 284, RUSSIA\\}

\date{July 1996}

\maketitle

\abstract{The tetrahedron equation
in a special substitution is reduced to a pair of 
pentagon and one ten-term equations. Various examples
of solutions are found. $O$-doubles of Novikov, 
which generalize the Heisenberg double of a Hopf 
algebra, provide a particular
algebraic solution to the problem.}
\vskip 2cm

\rightline{{\small E}N{\large S}{\Large L}{\large A}P{\small P}-L-611/96}

\end{titlepage}

\section{Introduction}

The Yang-Baxter equation (YBE) \cite{Yan,Bax} can be considered as
a tool for both constructing and solving integrable two-dimensional
models of statistical mechanics and quantum field theory \cite{Baxb,Fad}. 
Recent progress in understanding of the algebraic structure, lying behind
the YBE, has led to the
theory of quasi-triangular Hopf algebras \cite{Dri1}.

The tetrahedron (or three-simplex) equation (TE) \cite{Zam}
 has been introduced as a 
three-dimensional generalization of the YBE. Before describing
it in an abstract algebraic form, first consider an associative
unital algebra $A$, and define an important notation to be used 
throughout the paper. Namely, for each
set of integers $\{i_1,i_2,\ldots,i_m\}\subset \{1,2,\dots,n\}$, for $m<n$,
define a mapping 
\[
\tau_{i_1i_2\ldots i_m}\colon A^{\otimes m}\rightarrow A^{\otimes n}
\]
in such a way that
\[
a\otimes b\otimes\ldots\otimes c\mapsto 1\otimes\ldots\otimes
 a\otimes
\ldots\otimes b\otimes\ldots\otimes
 c\otimes\ldots\otimes1,
\]
where $a,b,\ldots,c$ in the r.h.s. stand on 
$i_1,i_2,\ldots, i_m$ -th positions, respectively, and unit
elements, on others. The 
notation to be used is as follows:
\be\label{not}
u_{i_1i_2\ldots i_m}=\tau_{i_1i_2\ldots i_m}(u),\quad u\in A^{\otimes m},
\ee
which means that the subscripts indicate the way
an element of the algebra $A^{\otimes m}$ is interpreted
as an element of the algebra $A^{\otimes n}$.

Next, we shall find it convenient to use the ``permutation operator'' 
$P$, which is an (additional) element in
$A^{\otimes2}$, defined by
\be\label{per}
Pa\otimes b=b\otimes a P,\quad P^2=1\otimes1, \quad a,b\in A.
\ee

The (constant) TE can be written as a 
nonlinear relation in $A^{\otimes6}$ on 
an invertible element $R\in A^{\otimes3}$ as follows:
\be\label{tet}
R_{123}R_{145}R_{246}R_{356}=R_{356}R_{246}R_{145}R_{123},
\ee
where we have used the notation (\ref{not}). One can introduce also
the higher-dimensional analogs of the YBE \cite{BS}, the $n$-simplex
equations. For example, the (constant) four-simplex equation (FSE)
is a relation in $A^{\otimes10}$
on an invertible element $B\in A^{\otimes4}$:
\be\label{4sim}
B_{0123}B_{0456}B_{1478}B_{2579}B_{3689}=
B_{3689}B_{2579}B_{1478}B_{0456}B_{0123}.
\ee
Many solutions have been found already for the TE, e.g.
\cite{Zam,BB1,BB2,KMS1,KMS2,Kor,Hie,MSS}, though, an adequate algebraic
framework (an analog of quasi-triangular Hopf algebras) is still missing.
Practically nothing is known for the higher-simplex equations.

The purpose of this paper is to make a step towards the algebraic theory
of the TE. Our main result is that one and the same system
of equations
\be\label{ss1}
S_{12}S_{13}S_{23}=S_{23}S_{12},
\ee
\be\label{ss2}
\overline S_{23}\overline S_{13}\overline S_{12}=
\overline S_{12}\overline S_{23},
\ee
in $A^{\otimes3}$ and 
\be\label{ss3}
\overline S_{12}S_{13}\overline S_{14}S_{24}\overline S_{34}=
S_{24}\overline S_{34}S_{14}\overline S_{12}S_{13},
\ee
in $A^{\otimes4}$ on elements $S,\overline S\in A^{\otimes2}$
implies both the TE for the combination
 \be\label{sans}
R^S_{123}=\overline S_{13} P_{23}S_{13},
\ee
and the FSE for the combination
\be\label{sans1}
B^S_{0123}=\overline S_{13}P_{01}P_{23}S_{13}.
\ee
We should warn, however, that formula (\ref{sans1}) is
too restrictive to give genuinely four-dimensional models. As a matter
of fact, it corresponds to non-interacting system of 
three-dimensional models. 

The manifest symmetry properties of equations 
(\ref{ss1}) -- (\ref{ss3}) are given by the following
transformations of $S$ and $\overline S$:
\be
S_{12}\leftrightarrow\overline S_{21},\quad 
S_{12}\leftrightarrow(\overline S_{12})^{-1}.
\ee
Equations (\ref{ss1}) and (\ref{ss2}) are the two forms of the 
celebrated pentagon equation (PE), which appears in various forms 
in representation theory of quantum groups as the Bidenharn-Elliot
identity for the 6j-symbols, in quantum conformal field theory as
an identity for the fusion matrices \cite{MS}, in quasi-Hopf algebras as
the consistency equation for the associator \cite{Dri2}. In the form
of (\ref{ss1}), (\ref{ss2}) the PE first appeared in the geometric
approach to three-dimensional integrable systems \cite{Mai1,Mai2}.
In \cite{Mai2} a reduction of the TE to the PE has been suggested. 
Finally, the
PE in the form (\ref{ss1}), (\ref{ss2}) was shown in \cite{Kas} to be
intimately related
with the Heisenberg double of a Hopf algebra \cite{RS,AF,Sem}.
In particular, using
the inclusion of the Drinfeld double into the tensor product of two
Heisenbergs, one can reduce the YBE to the PE.

Equation (\ref{ss3}), the ``ten-term'' relation, on
two different solutions of the PE, appears to be satisfied
by canonical elements in the $O$-double of a Hopf algebra,
introduced in \cite{Nov} as a generalization of the Heisenberg
double. Thus, the $O$-double is, probably, the simplest algebraic
framework for the TE.

The paper is organized as follows. In Sec.~\ref{sec1} particular 
solutions
for the PE, which generalize the solutions, associated with the
Heisenberg doubles of group algebras, are considered. The results of
this section are used in the next one for construction of particular 
solutions (typically
infinite dimensional) for the system (\ref{ss1}) -- (\ref{ss3}). 
In Sec.~\ref{sec2} the latter system is derived from the TE, and FSE, 
and the examples of solutions are described. In Sec.~\ref{sec3}
the $O$-double construction of a special class of solutions
is presented.

\section{Pentagon relation for a rational transformation}
\label{sec1}

As is shown in \cite{Kas}, for a group
$G$ the operator
\be\label{gr}
(S\varphi)(x,y)=\varphi(xy,y),\quad \varphi\in {\cal F}(G\times G),
\quad x,y\in G,
\ee 
in the space of functions on $G\times G$ satisfies PE (\ref{ss1}). 
Actually,
this is the ``coordinate'' representation for the canonical element
in the Heisenberg double of the group algebra.
In this section we generalize this result. Namely,
let $M$ be some set. Define an operator $S$ in the space
${\cal F}(M\times M)$:
\be\label{m0}
(S\varphi)(x,y)=\varphi(x\cdot y,x* y),\quad \varphi\in {\cal F}(M\times M),
\ee
for some mappings
\be\label{map}
M\times M\stackrel{\cdot }{\rightarrow} M,\quad
M\times M\stackrel{*}{\rightarrow} M.\quad
\ee
Let us call these the dot- and star-mapping, respectively.
Imposing now PE (\ref{ss1}) on $S$, we obtain the following equations:
\be\label{m1}
(x\cdot y)\cdot z=x\cdot (y\cdot z),
\ee
\be\label{m2}
(x* y)\cdot ((x\cdot y)* z)=x*(y\cdot z),
\ee
\be\label{m3}
(x* y)*((x\cdot y)* z)=y* z.
\ee
 If $M$ is a group with respect 
to the dot-mapping, then equation (\ref{m2}) implies
that
\be\label{star}
x*y=[x]^{-1}\cdot[x\cdot y],
\ee
where an invertible mapping $[\ ]\colon M\rightarrow M$ is given by 
\[
[x]=1*x.
\]
Substituting
(\ref{star}) into (\ref{m3}), we get the following functional equation
\be\label{fun}
\left[[x]^{-1}\cdot[x\cdot y]\right]=\lambda(x)\cdot[y],
\ee
with some function $\lambda\colon M\rightarrow M$. 
Putting $y=1$ in (\ref{fun}), we
obtain $\lambda(x)=1$, and therefore $[x]=x$. Thus, we have proved the
following proposition.
\begin{prop}
Let the set $M$ be a group with respect to the dot-mapping. 
Then, the star-mapping of the form $x*y=y$ is the only
solution to the system (\ref{m2})--(\ref{m3}).
\end{prop}
Remind, that the corresponding $S$-operator, given by (\ref{gr}), 
is associated with the Heisenberg double of the group algebra.
There are, however, other solutions, if the mappings under consideration
are defined only for elements in ``general position'', e.g.
birational mappings.
\begin{prop}
Let the set $M$ be a subset of an associative ring with unit,
and let the definitions
\be
x\cdot y=xy,\quad x*y=(1-x^\epsilon)^{-\epsilon}(1-(xy)^\epsilon)^\epsilon
,\quad \epsilon=\pm1,
\ee
 make sense in $M$. 
Then, equations (\ref{m1}) -- (\ref{m3}) are satisfied.
\end{prop}
The proof is straightforward.
\begin{prop}
Let $M=(0,1)\subset{\Bbb R}$ be the open unit interval with the dot-mapping
given by the multiplication in $\Bbb R$, and let the star-mapping
be continuously differentiable. Then, the system (\ref{m2})--(\ref{m3}) is  
satisfied iff
\be\label{dil}
x*y=y\left(\frac{1-x^{1/\alpha}}{1-(xy)^{1/\alpha}}\right)^\alpha,
\ee
where real $\alpha\ge 0$, and the case $\alpha=0$ is understood as a limit 
$\alpha\to 0^+$.
\end{prop}
{\it Proof}. It is straightforward to check that formula (\ref{dil}) 
solves the system (\ref{m2})--(\ref{m3}. Let us prove that it is the only 
continuously differentiable solution.

 Formulae (\ref{star}) and (\ref{fun}) are still valid, with
the mapping 
\[
[\ ]\colon M\to \Bbb R_+
\]
(rather than $M\to M$) being strictly 
increasing and continuously differentiable. Let us differentiate
(\ref{fun}) with respect to $y$. The result can be written as
\be\label{xx3}
w\left([xy]/[x]\right)=w(y)/w(xy),
\ee
where
\be\label{w}
w(x)=d\log x/d\log[x],
\ee
while differentiation of (\ref{fun}) with respect to $x$ together with
(\ref{xx3}) gives
\be\label{xx1}
w(xy)=w(x)-w(y)w(x)d\log\lambda(x)/d\log x.
\ee
Consistency of the last equation under the permutation 
$x\leftrightarrow y$ fixes
the derivative of the $\lambda$-function up to a real constant $c$:
\[
d\log\lambda(x)/d\log x=c-1/w(x).
\]
Plugging this back into (\ref{xx1}), we obtain the closed functional
equation
\be\label{xx2}
w(xy)=w(x)+w(y)-cw(x)w(y).
\ee
Here $c\ne0$ (otherwise $w(x)\sim\log x$, and eqs. (\ref{xx3}) and (\ref{w}) 
are in contradiction), therefore, equation (\ref{xx2}) can be rewritten in the
form
\[
1-cw(xy)=(1-cw(x))(1-cw(y)),
\]
the general continuous solution of which is well known:
\[
1-cw(x)=x^{1/\alpha},\quad 1/\alpha\in\Bbb R.
\]
Compatibility of eqs. (\ref{xx3}) and (\ref{w}) fixes
$c=1$. Thus,
\[
w(x)=1-x^{1/\alpha}>0 \Longrightarrow \alpha>0,
\]
where the first inequality follows from the definition (\ref{w}) of $w(x)$
and the fact that $[x]$ is strictly monotonically increasing function. 
Finally, solving the differential equation (\ref{w}), 
we complete the proof.

In conclusion, note that any two non-zero 
parameters $\alpha, \alpha'\ne0$ in (\ref{dil}) give
equivalent $S$-operators (\ref{m0}), consequently, we have only 
two inequivalent solutions for $M=(0,1)$, corresponding to $\alpha=0$,
and $\alpha=1$.

\section{Pentagon, ten-term, three-, and four-simplex relations}
\label{sec2}

Consider the following "ansatz" for a solution to equation (\ref{tet}):
\be\label{tans}
R^T_{123}=\overline T_{13} P_{23}T_{13},
\ee
for invertible elments $T,\ \overline T\in A^{\otimes2}$, with $P$ being
defined in (\ref{per}). 
\begin{prop}\label{10ter}
The TE (\ref{tet}) for an element $R$ of the form (\ref{tans}) 
is equivalent
to the existence of invertible elements $S,\ \overline S\in A^{\otimes2}$
such that the following equations are satisfied:
\be\label{st1}
S_{12}T_{13}T_{23}=T_{23}T_{12},\quad
\overline T_{23}\overline T_{13}\overline S_{12}=
\overline T_{12}\overline T_{23},
\ee
\be\label{st2}
\overline S_{12}T_{13}\overline T_{14}T_{24}\overline T_{34}=
T_{24}\overline T_{34}T_{14}\overline T_{12}S_{13}.
\ee
\end{prop}
{\it Proof}. Substituting (\ref{tans}) into (\ref{tet}), moving
all $P$-elements to the right, one can remove all of them from
the both sides of the equality simultaneously. 
The resulting identity can be rewritten in the form
\[
(\overline T_{16}^{-1}\overline T_{36}^{-1}\overline T_{13}\overline T_{36})
T_{12}\overline T_{15}T_{35}\overline T_{25}=T_{35}\overline T_{25}
T_{15}\overline T_{13}(T_{24}T_{12}T_{24}^{-1}T_{14}^{-1}).
\]
Here nontrivial elements in the subspaces 4 and 6 are contained only in the
expressions enclosed in brackets in the r.h.s. and l.h.s., respectively.
Consequently, these expressions should be trivial in the 
subspaces 4 and 6. In this way, we immediately 
come to the statements of the proposition.

\begin{prop}
Equations (\ref{st1}) and (\ref{st2}) imply that the 
elements $S$ and $\overline S$ satisfy equations 
(\ref{ss1}) --(\ref{ss3}).
\end{prop}
{\it Proof}. Relation (\ref{ss1}) follows from the identity
\[
S_{12}S_{13}S_{23}T_{14}T_{24}T_{34}=S_{23}S_{12}T_{14}T_{24}T_{34},
\]
which is proved by successive applications of the 
first identity from (\ref{st1}).
Relation (\ref{ss2}) is proved similarly. 
As for (\ref{ss3}), it
is a consequence of the identity
\[
\overline S_{12}S_{13}\overline S_{14}S_{24}\overline S_{34}
T_{15}T_{35}\overline T_{16}T_{26}
\overline T_{36}T_{46}\overline T_{56}=
\]
\[
=S_{24}\overline S_{34}S_{14}\overline S_{12}S_{13}
T_{15}T_{35}\overline T_{16}T_{26}\overline T_{36}T_{46}\overline T_{56},
\]
which can also  be proved by successive applications of the first identity
from (\ref{st1}) and (\ref{st2}).

The next proposition concerns the similar statement for the FSE.
\begin{prop}
The FSE (\ref{4sim}) for an element $B$ of the form (\ref{sans1})
is equivalent to equations (\ref{ss1}) -- (\ref{ss3}).
\end{prop}
The proof is similar to that of Proposition~\ref{10ter}.

Thus, relations (\ref{ss1}) -- (\ref{ss3}) enable us to
construct a special class of solutions for the TE and FSE. 
Let us turn to particular examples.
\begin{exam}
\rm
Take for the algebra $A$ the space of birational transformations 
${\rm Aut}(\Bbb C(x))$ of the field of rational expressions in one 
indeterminate $\Bbb C(x)$, and identify $A\otimes A\otimes\ldots$ with
 ${\rm Aut}(\Bbb C(x,y,\ldots))$. Interpreting ${\Bbb C}(x,y,\ldots)$ as 
the space of rational functions on ${\Bbb C}\times{\Bbb C}\times\ldots$,
define
\[
 (T\varphi)(x,y)=\varphi(xy,y-xy),\quad
\overline T=T^{-1}.
\]
Then, relations (\ref{st1}), (\ref{st2}) as well as (\ref{ss1}) -- (\ref{ss3})
are satisfied with
\[
 (S\varphi)(x,y)=\varphi(xy,[xy]/[x]),\quad[x]=x/(1-x),\quad
\overline S=S^{-1}.
\]
One can show that the corresponding element $R^T$ is equivalent to the 
solution $\Phi_0$ for the TE from \cite{Kas1}, which in turn is associated
with the three-dimensional Hirota equation of the discrete Toda system
\cite{Hir}.
\end{exam}
\begin{exam}
\rm
Let now  $\vec x=(x_1,x_2)$ be a pair of indeterminates, and put
$A={\rm Aut}(\Bbb C(\vec x))$, $A\otimes A\otimes\ldots$
being identified with ${\rm Aut}(\Bbb C(\vec x,\vec y,\ldots))$. Consider
the following rational mappings:
\[
\vec x\cdot\vec y=(x_1y_1,x_1y_2+x_2),\quad
\vec x*\vec y=[\vec x]^{-1}\cdot[\vec x\cdot\vec y],
\]
where
\[
\vec x^{-1}=(1/x_1,-x_2/x_1),\quad [\vec x]=(x_1/(1-x_1),x_2/(1-x_1)),
\]
and define
\[
(S\varphi)(\vec x,\vec y)=\varphi(\vec x\cdot\vec y,\vec x*\vec y),\quad
\overline S=S^{-1}.
\]
Then, equations (\ref{ss1}) -- (\ref{ss3}) are satisfied.
\end{exam}
\begin{exam}
\rm
Let algebra $A$ be the Heisenberg algebra, generated by elements
$\{ H,\Lambda,1\}$, satisfying the Heisenberg commutation relation
\[
\Lambda H-H\Lambda=1/h,
\]
with $h$ being a complex parameter with a positive real part.
Define the function 
\[
(x;q)_\infty=\prod_{n=0}^\infty(1-xq^n),
\]
where $q=\exp(-h)$, and put
\[
S=q^{H\otimes\Lambda}(-q^\Lambda\otimes q^{-H}q^{-\Lambda};q)_\infty^{-1}.
\]
Note, that this formula is a specialization of the canonical 
element in the Heisenberg double of the Borel subalgebra of 
$U_q(sl(2))$ quantum group,
see \cite{Kas}.
Now, both choices of $\overline S$, either
\[
\overline S=S^{-1},
\]
or
\[
\overline S= q^{-H\otimes\Lambda},
\]
solve the system (\ref{ss1}) -- (\ref{ss3}). The corresponding solutions 
(\ref{sans}) and (\ref{sans1}) to the TE and FSE first have been found in 
\cite{SBBMS}.
\end{exam}

\section{ $O$-double construction}
\label{sec3}

One particular class of solutions to the system (\ref{ss1}) -- (\ref{ss3})
is connected with $O$-doubles \cite{Nov}, 
which generalize the Heisenberg double of 
a Hopf algebra \cite{RS,AF,Sem}.

Consider elements $S$ and $\overline S$, satisfying (\ref{ss1}), (\ref{ss2})
but instead of (\ref{ss3}), impose  more restrictive set
of equations:
\[
S_{13}\overline S_{23}=\overline S_{23}S_{13},
\]
\be\label{co}
S_{12}\overline S_{13}\overline S_{23}=\overline S_{23}S_{12},\quad
S_{23}S_{13}\overline S_{12}=\overline S_{12}S_{23},
\ee
which imply also (\ref{ss3}). It appears that there is a general algebraic
structure, underlying equations (\ref{ss1}), (\ref{ss2}), and (\ref{co}).

Let $X$ be a Hopf algebra. In a linear basis $\{e_i\}$ 
the product, the  coproduct, the unit, the counit, and the antipode
take the form
\[
e_ie_j=m^k_{ij}e_k,\quad \Delta(e_i)=\mu_i^{jk}e_j\otimes e_k,
\]
\be\label{bas}
1=\varepsilon^ie_i,\quad\varepsilon(e_i)=\varepsilon_i,\quad
\gamma(e_i)=\gamma_i^je_j,
\ee
where summation over repeated indices is implied. Here $m^k_{ij}$,
$\mu_i^{jk}$, $\varepsilon^i$, $\varepsilon_i$, and $\gamma_i^j$ 
are numerical structure constants of the algebra.

Let $X^*$ be the dual Hopf algebra.
Following \cite{Nov}, consider an algebra $X^*XX^*$ ($O$-double), 
generated by right derivations
$R^*_x$, $x\in X$:
\be\label{al1}
R^*_x\colon X^*\rightarrow X^*,\quad \langle R^*_x(f),y\rangle=
\langle f,R_x(y)\rangle=\langle f,yx\rangle,
\ee
left, $L_f$, and right, $R_{\gamma^{-1}(g)}$, multiplications , $f,g\in X^*$:
\be\label{al2}
L_f,\ R_g\colon X^*\rightarrow X^*,\quad L_f(g)=R_g(f)=fg.
\ee

\begin{prop}
The algebra $X^*XX^*$ is an associative algebra, generated by elements
$\{e^i,e_j,\tilde e^k\}$, and the following defining relations:
\[
e^ie^j=\mu_k^{ij}e^k,\quad e_ie_j=m^k_{ij}e_k,\quad 
\tilde e^i\tilde e^j=\mu_k^{ij}\tilde e^k,
\]
\be\label{dub}
e_ie^j=m^j_{kl}\mu_i^{lm}e^ke_m,\quad 
\tilde e^ie_j=\mu_j^{kl}m^i_{lm}e_k\tilde e^m,\quad
e^i\tilde e^j=\tilde e^je^i.
\ee
\end{prop}
{\it Proof}. One has just to write the compositions of the operations, 
defined
in (\ref{al1}) and (\ref{al2}), for elements of the linear basis, 
using
(\ref{bas}) and the corresponding relations for the dual algebra.

\begin{prop}
Two canonical elements $S=e_i\otimes e^i$, 
$\overline S=e_i\otimes\tilde e^i$
in the $O$-double $X^*XX^*$ satisfy equations (\ref{ss1}), 
(\ref{ss2}), and (\ref{co}).
\end{prop}
The proof is straightforward through substitution of the canonical elements
into the relations to be proved, and application of formulae (\ref{dub}).

Thus , we have obtained a particular class of general algebraic solutions
to the system (\ref{ss1}) -- (\ref{ss3}), which in turn imply the 
TE for the element (\ref{sans}) and the FSE for the element (\ref{sans1}).

\section{Summary}

Solutions for the system of equations (\ref{ss1}) -- (\ref{ss3}) 
provide us both
with solutions for the three- and four-simplex equations (\ref{tet}) and 
(\ref{4sim}) through formulae (\ref{sans}) and (\ref{sans1}), respectively.

The $O$-double of a Hopf algebra provides an algebraic structure, underlying
the system (\ref{ss1}), (\ref{ss2}) and (\ref{co}), which implies
also (\ref{ss3}). Nevertheless, examples of solutions to the system
(\ref{ss1}) -- (\ref{ss3}), described in Sec.~\ref{sec2}, do not come
from the $O$-double construction. This suggests, that the latter
is a particular case of a more general algebraic structure, lying
behind the system (\ref{ss1}) -- (\ref{ss3}) itself.

\section{Acknowledgments}

The authors are indebted to A.Yu. Volkov for reading the manuscript and
helpful suggestions. It is a pleasure to thank also 
V.O. Tarasov, Yu.G. Stroganov, J.M. Maillet, L. Freidel, H.E. Boos for 
discussions.


\begin{thebibliography}{99}

\bibitem{AF}
  Alekseev, A.Yu.,  Faddeev, L.D.: $(T^*G)_t:$ a toy model for
 conformal field theory. 
Commun. Math. Phys. {\bf141}, 413--422 (1991) 

\bibitem{Baxb}
 Baxter, R.J.: Exactly solved models in statistical mechanics.
London: Academic Press 1982

\bibitem{Bax}
 Baxter, R.J.: Partition function of the eight-vertex lattice model.
Ann. Phys. {\bf 70}, 193--228 (1972) 

\bibitem{BB1}
 Bazhanov, V.V., Baxter, R.J.: 
New solvable lattice models in three dimensions.
J. Stat. Phys. {\bf 69}, 453--485 (1992) 

\bibitem{BB2}
Bazhanov, V.V.,  Baxter, R.J.:
Star-triangle relation for 
a three -~dimensional model.
J. Stat. Phys. {\bf71}, 839--864 (1993)

\bibitem{BS}
Bazhanov, V.V.,  Stroganov, Yu.G.:
Commutativity conditions for transfer 
matrices on a multidimensional lattice.
Theor. Mat. Fiz. {\bf 52}, 105--113 (1982) 
[English transl.: Theor. and Math. Phys. {\bf52}, 685--691 (1983)]

\bibitem{Dri1}
 Drinfeld, V.G.: Quantum groups. In: Proc. Int. Cong. Math., Berkeley 1987,
 P. 798--820

\bibitem{Dri2}
 Drinfeld, V.G.: Quasi-Hopf algebras. 
Algebra and Analysis {\bf1}, 114 (1989)  

\bibitem{Fad}
 Faddeev, L.D.: Quantum completely integrable models in field theory.
Sov. Sci. Rev. C{\bf1}, 107--155 (1980) 

\bibitem{Hie}
 Hietarinta, J.: Labelling schemes for tetrahedron equations 
and dualities between them.
 J. Phys. A{\bf27}, 5727--5748 (1994) 

\bibitem{Hir}
 Hirota, R.: Discrete analogue of a generalized Toda equation.
J. Phys. Soc. Jpn., {\bf50}, 3785--3791 (1981) 

\bibitem{Kas}
 Kashaev, R.M.: Heisenberg double and pentagon relation.
Preprint ENSLAPP-L-512/95, q-alg/9503005, 
to appear in Algebra and Analysis

\bibitem{Kas1}
 Kashaev, R.M.: On discrete 3-dimensional equations 
associated with the local Yang-Baxter relation.
Preprint ENSLAPP-L-569/95, solv-int/9512005, 
to appear in Lett. Math. Phys.

\bibitem{KMS1}
 Kashaev, R.M., Mangazeev, V.V., Stroganov, Yu.G.: Spatial symmetry, 
local integrability, and tetrahedron equation 
in the Baxter-Bazhanov model.
Int. J. Mod. Phys. A{\bf8}, 587 (1993) 

\bibitem{KMS2}
Kashaev, R.M., Mangazeev, V.V., Stroganov, Yu.G.:
Star-square and tetrahedron equations in the 
Baxter-Bazhanov model. 
Int. J. Mod. Phys. A{\bf8}, 1399 (1993) 

\bibitem{Kor}
 Korepanov, I.G.: Tetrahedral Zamolodchikov algebras 
corresponding to Baxter's $L$-operator.
Commun. Math. Phys. {\bf154}, 85--97 (1993)

\bibitem{Mai1}
 Maillet, J.M.: Integrable systems and gauge theories.
 Nucl. Phys. (Proc. Suppl.) B{\bf18}, 212--241 (1990) 

\bibitem{Mai2}
 Maillet, J.M.: On Pentagon and Tetrahedron equations. 
Algebra and Analysis {\bf6}, 206 (1994)

\bibitem{MSS}
Mangazeev, V.V., Sergeev, S.M., Stroganov, Yu.G.:
New solutions of vertex type tetrahedron equations.
Mod. Phys. Lett. A{\bf10}, 279--287 (1995) 

\bibitem{MS}
Moore, G., Seiberg, N.: Classical and quantum conformal field theory.
Commun. Math. Phys. {\bf123}, 177--255 (1989) 

\bibitem{Nov}
 Novikov, S.P.: Various doublings of Hopf algebras.
 Algebras of operators on quantum groups. Complex cobordisms.
 Usp. Math. Nauk {\bf47},  No. 5(287), 189-190 (1992), 
transl. in Russian Math. Surveys {\bf47}, No. 5, 198--199 (1992)

\bibitem{RS}
 Reshetikhin, N.Yu., Semenov-Tian-Shansky, M.A.:
Central extensions of quantum current groups.
 Lett. Math. Phys. {\bf19}, 133--142 (1990) 

\bibitem{Sem}
  Semenov-Tian-Shansky, M.A.: Poisson-Lie qroups. 
The quantum duality principle and the twisted quantum double. 
Theor. Math. Phys.{\bf93}, 302--329 (1992), 
transl. in: Theoret. and Math. Phys. {\bf93}, No. 2, 1292--1307
(1992) 

\bibitem{SBBMS}
 Sergeev, S.M., Bazhanov, V.V., Boos, H.E., 
Mangazeev, V.V., Stroganov, Yu.G.: 
Quantum dilogarithm and tetrahedron equation.
Preprint IHEP 95-129

\bibitem{Yan}
 Yang, C.N.: Some exact results for the many-body 
problem in one dimension with 
repulsive delta-function interaction. 
Phys. Rev. Lett. {\bf 19}, 1312--1314 (1967)

\bibitem{Zam}
 Zamolodchikov, A.B.:
Tetrahedron equations and the relativistic 
$S$-matrix of straight-strings
in $2+1$ dimensions.
Commun. Math. Phys. {\bf79}, 489--505 (1981)  
 

\end{thebibliography}
\end{document}